\documentclass[11pt]{article}
\usepackage[margin=1.5in]{geometry} 

\newtheorem{theorem}{Theorem}[section]
\newtheorem{pro}{Proposition}[section]

\newtheorem{definition}{Definition}[section]
\newtheorem{remark}{Remark}[section]
\newtheorem{cor}{Corollary}[section]

\newcommand{\proof}[1]{\noindent{\it\bf Proof:#1\ }}
\newcommand{\QED}{QED}

\begin{document}

\title{   On the Action of   Reparametrization Group  on the Space of  $L_k^p$-maps I}
\author{Preliminary Version\\
\\
Gang Liu }
\date{December 3,  2013}
\maketitle

\begin{large}
\section{Introduction}
  
Let $(M, g)$ be a  compact $C^{\infty}$-manifold  with a Riemaniann metric
  $g$ and    $\Sigma\simeq  S^2$ with  its standard complex structure and round metric.  The group $G_0\simeq PSL (2, {\bf C})$ acts on  $\Sigma$ as  the reparametrization group. Let $G_i, i=1,2,$ be the subgroup of $G_0$ preserving one or two fixed marking points $x_1$ and $x_2$ of $\Sigma$. We will use $G$ to denote any of these three groups if there is no confusion.
Fix $k$ and $p$ such that $m_0=k-\frac{2}{p}>1$. Let  ${ {\cal M}}={{\cal M}}_{k,p}(\Sigma, M)$ be the space of 
 $L_k^p$-maps from $\Sigma$ to $M$.
Denote  the subspace of non-trivial 
 $L_k^p$-maps   by ${ {\cal M}^*}={{\cal M}}^*_{k,p}(\Sigma, M)$. The  reparametrization group 
$G$  acts on  ${\cal M}$. The purpose of this paper is  to prove the following theorem.

\begin{theorem}

  The  action  of $G$  on  ${ {\cal M}}^*$ is  proper in the following sense: 
  given any two maps $f_1$ and $f_2$ in ${ {\cal M}}^*$, there exist small neighbourhoods $U_{\epsilon_1}(f_1)$ and $U_{\epsilon_2}(f_2)$ of $f_1$ and $f_2$ in ${ {\cal M}}^*$ and two compact sets $K_{f_1}$ and $K_{f_2}$ in $G$ such that
  
  \vspace{2mm}
\noindent
  (a) for any $h$ in $U_{\epsilon_1}(f_1)$ and $g$ in $G\setminus K_{f_1}$, 
   $h\circ g$ is not in $U_{\epsilon_2}(f_2)$;
  
   \vspace{2mm}
\noindent
(b) the corresponding statement holds for any $h$ in $U_{\epsilon_2}(f_2)$ and $g$ in $G\setminus K_{f_2}$.
  
  \end{theorem}

\vspace{2mm}
\noindent  
 By taking $f=f_1=f_2$, we get the following theorem.

\begin{theorem}

For any non-constant $L_k^p$-map $f$, there exist a small neighbourhood $U_{\epsilon}(f)$ in ${ {\cal M}}^*$ and a  compact set $K_{f}$  in $G$ such that for any $h$ in $U_{\epsilon}(f)$ and $g$ in $G\setminus K_{f}$, 
   $h\circ g$ is not in $U_{\epsilon}(f)$. 
  
  In particular, for any such $g$
  in $G\setminus K_{f}$,  $f\circ g$ is not in $U_{\epsilon}(f)$.
In other words, the obit of $f$, $O_f=G\cdot f$, can not come back to the sufficient small neighbourhood $U_{\epsilon}(f)$ forever even $G$ is non-compact.

\end{theorem}

A corollary of this is the following.

 \begin{cor} 
 For any non-constant $L_k^p$-map $f$,  its stabilizer $stab_f$ is always a compact
  subgroup of $G$.
 \end{cor}

Further consequences of the Theorem 1.1 will be discussed in Sec. 3.

  The theorems stated above seem unreasonable. We now  state  three  simple  facts  and an immediate consequence of theirs which is even more counter-intuitive.

\vspace{2mm}
\noindent  
(A)   The  $L^2$-energy $E(f)=\int_{\Sigma}|df|^2\cdot dvol_{\Sigma}$ is  conformally invariant. That is $E(f\circ g)=E(f)$ for any $g\in G.$

 \vspace{2mm}
\noindent 
\vspace{2mm}
\noindent 
\vspace{2mm}
\noindent 
(B)  For two $L_k^p$-maps $f$ and $h$, the difference of their energies
$|E(f)-E(h)|$ is bounded by $C(\|f\|_{k, p}+\|h\|_{k, p})\|f-h\|_{k, p} $ for some constant $C$.  These two properties imply the following well-known fact.

\vspace{2mm}
\noindent 
\vspace{2mm}
\noindent 
\vspace{2mm}
\noindent 
$(C)$ The energy function $E:{\cal M}\rightarrow {\bf R}$ is continuous and $G$-invariant.

A  consequence of these facts is the following proposition.

\vspace{2mm}
\noindent

\begin{pro}
Given any two $L_k^p$-maps $f_1$ and $f_2$ with $E(f_1)\not =E(f_2)$, there exit   $G$-neighbourhoods  $  {\bf W}(f_1)$ of $f_1$ and $ {\bf W}(f_2)$ of $f_2$  which  do  not intersect.

 In particular if $f_1$ is a constant map and $f_2$ is not, then $E(f_1)=0\not =E(f_2)$ and the above conclusion  holds.

\end{pro}

\proof

Note that the condition $E(f_1)\not =E(f_2)$ implies that $f_1$ and $f_2$ are not in the same $G$-orbit. 
We may assume that $E(f_1)<E(f_2)$. For any $E(f_1)<c<E(f_2), $ since the energy function $E:{\cal M}\rightarrow {\bf R}$ is continuous and $G$-invariant, the inverse images $E^{-1}((-\infty, c))$ and $E^{-1}((c, \infty))$, denoted by  ${\bf W}(f_1)$ and $ {\bf W}_2(f_2)$,  are two open
$G$-sets in ${\cal M}$ containing $f_1$ and $f_2$ respectively. Clearly ${\bf W}(f_1)$ and $ {\bf W}_2(f_2)$ do not intersect.

\QED

 A corollary of this is the following Proposition.

\begin{pro}
Given any two $L_k^p$-maps $f_1$ and $f_2$ with $E(f_1)\not =E(f_2)$, there exit  two $G$-neighbourhoods  $ G\cdot  U_{\epsilon_1}(f_1)$ and $ G\cdot  U_{\epsilon_2}(f_2)$ which  do  not intersect.

 In particular if $f_1$ is a constant map and $f_2$ is not, the conclusion still holds.

\end{pro}

\proof

\vspace{2mm}
\noindent 
\vspace{2mm}
\noindent 
Proof I:

Assume the above Proposition is true. Let $U_{\epsilon_i}(f_i), i=1,2 $  be an open neighbourhood of $f_i$ contained in  $  {\bf W}(f_i).$
Since  $  {\bf W}(f_i)$ is an open $G$-set,  $G\cdot U_{\epsilon_i}(f_i)$ is still contained in  $  {\bf W}(f_i).$

\QED

\vspace{2mm}
\noindent 
\vspace{2mm}
\noindent 
Proof II:

We now give a direct and more computational proof without using the above Proposition.

For completeness, we give a detail proof for (B) above first.

For two $L_k^p$-maps $f$ and $h$,  $$|E(f)-E(h)|=|\Sigma_{i, j}\int_{\Sigma}({\partial}_if^j)^2-({\partial}_ih^j)^2\cdot d{\bf x}|$$ $$\leq
\Sigma_{i, j}\int_{\Sigma} (|{\partial}_if^j)|+|{\partial}_ih^j|)\cdot 
|{\partial}_if^j-{\partial}_ih^j|d{\bf x}$$ $$\leq C(\|f\|_{C^1}+\|h\|_{C^1})\|f-h\|_{C^1} $$ $$ \leq C\cdot (\|f\|_{k, p}+\|h\|_{ k, p})\|f-h\|_{k, p}  $$ $$  \leq C\cdot (2\|f\|_{k, p}+\|f-h\|_{ k, p})\|f-h\|_{k, p} .$$

\vspace{2mm}
\noindent 
Here  ${\bf x}=(x_1, x_2)$ is a  conformal coordinate chart  for $\Sigma$ and metric on $M$  is the one induced from an embedding of $M$
into some Euclidean  space.

\noindent 
 This  proves (B).
 
Now   let $f_1$  and $f_2$ be two $L_k^p$-maps such that $E(f_1)\not= E(f_2).$ 
Then there are neighbourhoods $U_{\epsilon_1}(f_1)$ and $U_{\epsilon_2}(f_2)$ such that $U_{\epsilon_1}(f_1)\cap U_{\epsilon_2}(f_2)$ is empty. 

 For any $k_i\in G\cdot  U_{\epsilon_i}(f_i), i=1,2$, write $k_i=h_i\circ g_i$ with $h_i\in  U_{\epsilon_i}(f_i)$ and $g_i\in G.$ 
 
 Then since the automorphism  $g_i$ of $\Sigma$ is conformal, by (A),  we have $$|E(k_i)-E(f_i)|=|E(h_i\circ g_i)-E(f_i)|$$ $$ =|E(h_i)-E(f_i)|\leq 
 C\cdot (2\|f_i\|_{k, p}+\|f_i-h_i\|_{ k, p})\|f_i-h_i\|_{k, p} $$
 $$\leq 
 C\cdot (2\|f_i\|_{k, p}+\epsilon_i) \epsilon_i.$$

Denote  $|E(f_1)- E(f_2)|$  by $\delta>0$ and $\epsilon_1+\epsilon_2$ by $\epsilon$.  Then $$|E(k_1)-E(k_2)|=|(E(k_1)-E(f_1))-(E(k_2)-E(f_2))+(E(f_1)-E(f_2))|$$
$$\geq |(E(f_1)-E(f_2))|-|(E(k_1)-E(f_1))|-|(E(k_2)-E(f_2))|$$

$$\geq  |(E(f_1)-E(f_2))|-2C
(\|f_1\|_{k, p}+\|f_2\|_{k, p}+\epsilon_1+\epsilon_2)( \epsilon_1+\epsilon_2)$$.

$$\geq \delta -2C\epsilon
(\|f_1\|_{k, p}+\|f_2\|_{k, p}+\epsilon). $$

\vspace{2mm}
\noindent 
if we choose $\epsilon$ small enough so that $\delta -2C\epsilon
(\|f_1\|_{k, p}+\|f_2\|_{k, p}+\epsilon)>0, $ then $|E(k_1)-E(k_2)|>0$ for any
$k_1 \in G\cdot  U_{\epsilon_1}(f_1)$ and $k_2 \in G\cdot  U_{\epsilon_2}(f_2)$. This implies that  $U_{\epsilon_1}(f_1)$ and $U_{\epsilon_2}(f_2)$ do not intersect each other.

\QED

The last statement of the propositions seems more against our intuition and experience in Gromov-Witten
theory.  Since the constant maps   have the whole non-compact reparametrization group $G$ as their stabilizer, they are certainly unstable in any reasonable sense.  The  common practice  in $GW$-theory, or  in any theory on how to form  "good" quotient spaces like GIT,  is to  exclude such unstable points even for much rigid situation like moduli space of $J$-holomorphic maps. Yet our simple argument above shows that in this particular case, these unstable  points are harmless for our purpose. In fact, in Sec. 3,  we will show that the $G$-space ${\cal M}$ of all $L_k^p$-maps is $G$-Hausdroff in the sense that any two points not in the same $G$-orbit  are always separable  by two $G$-neighbourhoods.
This makes the situation even worse. To make the discussion here more compatible with the general belief in Gromov-Witten theory, we note that allowing the appearances of the constant  maps as  trivial bubbles does cause trouble and produces non-Hausdroff quotient  spaces in GW-theory. But this only occurs when there are changes of  topological types of the domains
$\Sigma$. As long as  the node curves $\Sigma$ stay  in the same stratum, we always have the above conclusion for $L_k^p$-maps.

In fact, in one of the sequels of this paper [L1], we will show that the in GW-theory, only way that "bad" thing can happen for the  space of $L_k^p$-maps is the trivial way that  a sequence of such maps "converges" to a $L_k^p$-map which has
extra trivial bubbles. A typical well-known example of such a sequence is $f_n=f\circ g_n:{\bf CP}^1\rightarrow M$ with $g_n\in {\bf C}^*\subset PSL(2, {\bf C})$ going to zero or infinity  of ${\bf C}^*$ so that the sequence produces
trivial bubbles in a trivial manner. Note that in this example, if we allow the   limits  of the sequence, the limits are not in the  space ${\cal M}.$

 This paper is organized a follows.
 
 Sec. 2 proves the main theorems and the Corollary 1.1. The proof here is not
 most effective  and arguments are often repeated.   It is written only  using 
 elementary arguments and insisting on all details. A more effective proof for more general situations is in [L0].
 
 Sec. 3 proves  more corollaries of the Theorem 1.1.
 
 Sec.4 states the generalizations of the main theorem to the following two cases:
 (i) the domain $\Sigma$ is a smooth projective manifold and the reparametrization group $G$ is a connected reductive group acting on  $\Sigma$ algebraically; (ii)  $\Sigma$ is a smooth and compact Kahler  manifold and the group $G$ is a connected reductive group acting on  $\Sigma$ holomorphically such that the action comes from a moment map. The proofs for theses generalizations will be given in [L0].
\vspace{2mm}
\noindent  

 Sec 5.  defines a pseudo-moment map for $K=SU(2, {\bf C})$ acting on ${\cal M}^*$ above, as well as a corresponding pseudo-moment map for the maximum compact subgroup $K$ in the reductive group $G$  in the above general cases (i) and (ii). The discussion here  is very brief.  Such a pseudo-moment map can be used to give a global slice for the non-compact directions  of the $G$-action on  ${\cal M}^*$ at least for the case of an open neighbourhood  of the moduli space of $J$-holomorphic maps in  GW-theory. The implications of the existence of such a global slicing for the problem of regularizing the moduli space in GW-theory  and the  possibility to have a  GW-theory  with $SU(2)$-action  will be discussed in the revised version of this paper.

\section{The Proof of the Main Theorem}
\vspace{2mm}
\noindent  
\vspace{2mm}
\noindent  
 $\bullet$ The proof of the corollary 1.1:

\vspace{2mm}
\noindent 
Denote the action map by $\Psi: G\times {\cal M}\rightarrow {\cal M}$.
We will assume the well-known fact that the orbit map $\Psi(-, f):G\rightarrow {\cal M}$  is continuous.
This implies that $Stab_f$ is closed in $G$.

\noindent 
$\bullet$ $\bullet$ Case (I): 

We start  with the case that $(\Sigma,x_1, x_2)\simeq ({\bf P}^1, 0, {\infty})=({\bf C}\cup \{{\infty}\}, 0, {\infty}).$ In this case, the group $G=G_2$ preserving
 $x_1$ and $x_2$ is $G=C^*=\{a \,|a\in {\bf C}, a \not = 0\}$. For any $a\in G$
 and $z\in {\bf C}\subset {\bf C}\cup \{{\infty}\}$, the action $\Psi ( a, z)=a\cdot z.$

Assume that $Stab_f\subset G$ is not compact, then there is a sequence $\{a_n\}_{n=0}^{\infty}$ in $Stab_f$ such that either $lim_{n\mapsto \infty} a_n=0$ or $lim_{n\mapsto \infty} a_n=\infty.$ Using the automorphism of $\Sigma$ given by $w=\frac {1}{z},$ we only need to consider one of the cases.

We assume that $lim_{n\mapsto \infty} a_n=0$. In fact, the other case is included
as a special case in our treatment for case (II) below.

Given a non-constant map $f:\Sigma\rightarrow M$, let $\delta=\delta_f>0$ be its diameter. Then there are two points $y_1$ and $y_2$ in $\Sigma$  such that  the distance $d(f(y_1), f(y_2))=\delta.$ In the case that one of them, $y_2$ for instance is equal to ${\infty}$, we may chose $y_2'\not = \infty$ such that  $d(f(y_1), f(y'_2))>\delta/2.$ Therefore,  there is a sufficient large disc $D_R$ of radius $R$ centred  at $x_1=0$ in ${\bf C}\subset \Sigma= {\bf C}\cup \{{\infty}\}$ such that the diameter of $f(D_R)>\delta/2.$ 
 
 On the other hand, the continuity of $f$ implies that there is a $\gamma >0 $ such that for any disc $D_r$ centred  at $0$  with  radius $r\leq\gamma , $  
 the diameter of $f(D_r)<\delta/2000.$ 
 
 Denote $a_n\in G$ above by $g_n:\Sigma\rightarrow \Sigma.$ Since $a_n\rightarrow
 0,$ for $n$ sufficiently large, we have that $g_n(D_R)$ is contained in $D_r$ with $r<\gamma.$

 Since $g_n$ is in $Stab_f$, we have $f=f\circ g_n$ for all $n$. Therefore, we have that for  $n$ sufficiently large,  $\delta/2<$
 the diameter of $f(D_R)= $ the diameter of $ f\circ g_n (D_R)\leq $
 the diameter of $ f (D_r)\leq \delta/2000.$ We get a contradiction. This finishes
  the proof of the corollary for case (I).
 
 \vspace{2mm}
\noindent  
\vspace{2mm}
\noindent  
 $\bullet$ $\bullet$ Case (II):

In this case, $G=G_1$, $ (\Sigma, x)=({\bf C}\cup \{{\infty}\}, \infty)$. 
Each element $g\in G$ has the form $g(z)=az+b$ with $a\not =0$ and $a, b$ in ${\bf C}$. We will write it as $g(z)=a(z-c).$ Note that $G_1$ contains all translations.
 
 Given a non-constant $L_k^p$-map $f, $ let $T=T_b$ defined by $T(z)=z+b$ be a translation
   in $G_1$, then the subgroups $Stab_f$ and $Stab_{f\circ T}$ are conjugate
   each other in $G$ by $T$. Therefore, by using a translation if necessary, we may assume that, $f(0)\not= f(\infty).$
   
 Now assume that $\{g_n\}_{n=1}^{\infty}$ given by $g_n(z)=a_n(z-c_n)$ be a sequence in $Stab_f$.  
 
 \vspace{2mm}
\noindent  
\vspace{2mm}
\noindent  
 $\bullet$ $\bullet$ $\bullet$ Claim:
 $|c_n|$ is bounded.
 
 \proof
 
 Assume that the claim is not true. Then there is a subsequence of $\{g_n\}_{n=1}^{\infty}$, denoted by the same notation, such that  $c_n\rightarrow \infty$  as $n\rightarrow \infty.$
 Since $g_n$ is in $Stab_f$, we have $f=f\circ g_n$. Therefore,
 $f(c_n)=f\circ g_n(c_n)=f(a_n(c_n-c_n))=f(0).$ The continuity of $f$ at $\infty$ implies that
 $f(\infty)=\lim_{n\mapsto\infty}f(c_n)=f(0),$ which contradicts to our assumption.
 
 \QED
 
 Therefore, by taking a subsequence, we may assume that $\lim_{n\mapsto\infty}c_n=c \in {\bf C}.$
 
 Now assume that $Stab_f$ is not compact. This implies that there is a sequence
 $\{g_n\}_{n=1}^{\infty}$ as  above given by $g_n(z)=a_n(z-c_n)$ which has no convergent subsequences in $Stab_f$. Since $Stab_f$ is closed in $G$, this can only happen if the sequence $\{a_n\}_{n=1}^{\infty}$ in $C^*$  has a subsequence  either going  to zero or to infinity. Now  we are almost in the same situation as the case (I) except there is a bounded shifting of the origin by $c_n.$
 
Like  the case (I), the arguments for  the two cases are similar. This time we  assume that $\lim_{n\mapsto\infty}a_n=\infty.$
 
 Still assume that  diameter of $ f =\delta$. Then there is a positive
 $\rho$, such that when $\rho$ is small enough, the complement in $\Sigma$ of the disc
 $D_{\rho}(c)$ of radius $\rho$ centred at $c$, denoted by $D^c_{\rho}(c)$,
  has the property that the diameter of $f(D^c_{\rho}(c))>\delta/2.$ 
 
  For any $z\in D_{\rho/2}(c_n), $  we have that for $n$ large enough, $$|z-c|\leq |z-c_n|+|c_n-c|
 < \rho/2+\rho/2=\rho.$$
 
 In other words, for large $n$, $D_{\rho/2}(c_n)\subset  D_{\rho}(c)$. Hence
 we have that   $D_{\rho}^c(c)\subset  D_{\rho/2}^c(c_n).$

 Now let ${\tilde D}_r(\infty)$ be a small disc in $\Sigma={\bf C}\cup \{\infty\}$ centred at $\infty$ of radius $r$ measured in the Fubini-Study
 metric such that the diameter of $ f ({\tilde D}_r(\infty))\leq \delta/2000.$

 Note that there is a large disc $D_R$ of radius $R$ centred  at $0$ in ${\bf C}$
 such that $ {\tilde D}_r(\infty)=D^c_R$.
 
 Now for any $z$ in $D_{\rho}^c(c)$, when $n$ is large enough, we have that
 $$|g_n(z)|=|a_n|\cdot |z-c_n|\geq |a_n|\cdot (|z-c|-|c_n-c|)$$ 
 $$\geq |a_n|\cdot (\rho-|c_n-c|)\geq |a_n|\cdot (\rho-\rho/2)\geq R.$$
 
 This implies that for such $n$, $g_n(D_{\rho}^c(c))\subset D^c_R={\tilde D}_r(\infty) .$ Therefore for such $n$, we have that 
 $\delta/2<$ the diameter of $f(D^c_{\rho}(c))=$ the diameter of $f\circ g_n((D^c_{\rho}(c))<$ the diameter of $f(D^c_R))=$ the diameter of $f({\tilde D}_r(\infty)) \leq \delta/2000.$

This finishes the proof of the case (II).
\QED 
 
Before going any further, we make a few remarks on the two case above.

\vspace{2mm}
\noindent  
\vspace{2mm}
\noindent  
 $\bullet$ Remark:
 
\vspace{2mm}
\noindent  
(i) The group $G=G_0=PSL(2, {\bf C})=SL(2, {\bf C})/\{\pm 1\}.$ 
 Since $SL(2, {\bf C})$ is reductive, the case (I) above is typical for $G_0$ in the sense that the proofs of  the Theorem 1.1 and above corollary for $G_0$ can be reduced to  the case $G_2=C^*.$ The similar situation happens  in GIT known  as the Hilbert-Mumford criterion. Indeed, for any non-scalar $g \in SL(2, {\bf C}),$ upto a
   ${\bf Z}_2$-action,
 we have an unique decomposition in  $SL(2, {\bf C}),$ $g=h\cdot  u$ with $u\in SU(2)$ and $h$ being self-adjoint. Therefore, we have that  $g=u_1\cdot D(a)\cdot u_2$. Here  $u_i, i=1,2,$ are  in $SU(2, {\bf C})$ and $D(a)$ is the diagonal matrix with entries $a$ and $a^{-1}.$ This decomposition is also essentially unique upto some obvious  $SU(2)$-actions. Clearly,under the projection $ SL(2, {\bf C})\rightarrow PSL(2, {\bf C})$,  $SU(2, {\bf C})$ becomes the double covering of  the $SO(3)$ and the collection of all  $D(a)$ above maps onto $C^*$ sending $D(a)$ to $a^2$ (or $a^{-2}$).

Therefore any non-compact  sequence $\{[g_n]\}_{n=1}^{\infty}$ in $G=G_0$ with $g_n\in $ $SL(2, {\bf C})$    has the form
$g_n=u_n\cdot D(a_n)\cdot v_n$ with   $u_n$ and $v_n$   in $SU(2, {\bf C})$ and $a^2_n$ in $C^*$ such that  $\{a_n\}_{n=1}^{\infty}$ is a non-compact sequence
 in $C^*$.  From our proof in case (II) above, we already have the general idea on how to deal with the compact sequences like  $\{u_n\}_{n=1}^{\infty}$ and  $\{v_n\}_{n=1}^{\infty}$ here. This essentially reduces the proof of the Theorem 1.1 and its corollary for the case $G=G_0$ to the case (I) with $G=C^*.$
 Therefore for the proof of the Theorem 1.1 and Theorem 1.2 below in this section, we will only deal with the two cases as above.
 The case for general reductive group including $G_0$ is treated in [L0].

 \vspace{2mm}
\noindent  
\vspace{2mm}
\noindent  
 
\noindent  
(ii) For the applications in Gromov-Witten  and Floer theories, it is sufficient
to only consider above two cases. The $PSL(2, {\bf C})$-action on the top stratum  of the moduli space of stable maps used in genus zero GW-invariants can be removed by putting  constrains on the  stable maps with at least three marked points.

 \vspace{2mm}
\noindent  
\vspace{2mm}
\noindent  
 
\noindent

The proof above for Corollary 1.1  with a few small modifications implies the proofs of the two main theorems. It is sufficient to prove Theorem 1.1 (a).

\vspace{2mm}
\noindent  
\vspace{2mm}
\noindent  
 $\bullet$ The proof of the Theorem 1.1(a):

\vspace{2mm}
\noindent 
\noindent  
 Assume that the Theorem 1.1 (a)  is not true. Then for any small neighbourhoods
 $U_{\epsilon_i}(f_i), i=1, 2 $ and any nested  sequences of compact sets $K_1\subset K_2\subset \cdots \subset K_n\cdots$ in $G$, there  are  sequences
 $\{g_n\}_{n=1}^{\infty}$ in $G$ and $\{h_n\}_{n=1}^{\infty}$ in $U_{\epsilon_1}(f_1)$ such that (a) $g_n$ is not in $K_n$; (b) $h_n\circ
 g_n$ is in  $U_{\epsilon_2}(f_2).$ Here $\epsilon_i, i=1, 2 $ and $K_n, n=1, \cdots$ will be decided later in the proof.
 
\vspace{2mm}
\noindent  
 $\bullet$ $\bullet$ Case (I):

 We have   $G=G_2=C^*$.  Choose $K_n\subset C^*$ to be $\{a\in C^*\,|\,\frac {1}{n} \leq |a|\leq n\}$. Then the condition (a) above implies that for $g_n$ in $G$
  with $g_n(z)=a_nz$,  either $lim_{n\mapsto \infty} a_n=0$ or $lim_{n\mapsto \infty} a_n=\infty.$ 
  
As before, we only consider the case that $lim_{n\mapsto \infty} a_n=0$.

We have  already proved that there is a sufficient large disc $D_R$ of radius $R$ centred  at $0$ in ${\bf C}\subset \Sigma= {\bf C}\cup \{{\infty}\}$ such that the diameter of $f_2(D_R)>\delta_2/2>0, $  where $\delta_2$ is the diameter of the image of the non-constant map $f_2$. In particular, there are two points $y_1$ and 
$y_2$ in $D_R$ such that the distance $d(f_2(y_1), f_2(y_2))>\delta_2/2.$ Since
$h_n\circ g_n $ is in $U_{\epsilon_2}(f_2), $ we have $d(f_2(y_i), h_n\circ g_n(y_i))\leq C'_1\|f_2- h_n\circ g_n\|_{C^0},\,  i=1, 2 $ for some constant $C'_1$. Note that to make sense of the expression $h_n\circ g_n-f_2$ here, 
we have used the standard  exponential coordinate on $U_{\epsilon_2}(f_2)$
so that $ h_n\circ g_n=Exp_{f_2}\xi_n$ with $\xi_n\in L_k^p(\Sigma,  f_2^*(TM))$.
Then $h_n\circ g_n-f_2$ is defined to be $\xi_n.$
 Another way  to deal with this  is to embed  $M$ into some ${\bf R}^m$ so that
  ${\cal M}_{k, p}$ is contained in the Banach space ${\cal M}_{k, p}(\Sigma, {\bf R}^m)$ of $L_k^p$-maps from $\Sigma$ to ${\bf R}^m$. 
  
By our assumption,  $C'_1\|f_2- h_n\circ g_n\|_{C^0}\leq C_1\|f_2- h\circ g_n\|_{k, p}\leq C_1 \epsilon_2, $ which implies that $d(f_2(y_i), h_n\circ g_n(y_i))\leq C_1\epsilon_2, i=1,2 .$  We have 

$$d(h_n\circ g_n (y_1), h_n\circ g_n (y_2))$$ $$ \geq  d(f_2(y_1),  f_2(y_2))-d(f_2(y_1),  h_n\circ g_n (y_1))-d( h_n\circ g_n (y_2)  , f_2(y_2))$$ $$ \geq \delta_2/2-2C_1\epsilon_2. $$  Therefore, we have  
 the diameter of $h_n\circ g_n (D_R)\geq \delta_2/2-2C_1\epsilon_2 $ for any $n$ provided that we choose $2C_1\epsilon_2<< \delta/2.$
 
 As before, the continuity of $f_1$ implies that there is a $\gamma >0 $ such that for any disc $D_r$ centred  at $0$  with  radius $r\leq\gamma , $  
 the diameter of $f_1(D_r)<\delta_1/N.$  Here $N$ is fixed depending on $\delta_i, \epsilon_i,  i=1,2.$ It will be determined later in the proof.
 
 We  have proved that  since $a_n\rightarrow
 0,$ for $n$ sufficiently large (depending on $N$ above), we have that $g_n(D_R)$ is contained in $D_r$ with $r<\gamma.$ Since $h_n$ is in $U_{\epsilon_1}(f_1), $ we have that when $n$ is large enough, for any $y$ in  $D_R$, $$d(f_1(0), h_n\circ g_n(y))\leq d(f_1(0), f_1\circ g_n(y))+d(f_1\circ g_n(y) , h_n\circ g_n(y))$$ $$ \leq \delta_1/N+C_2 \|f_1- h_n\|_{C^0}\leq \delta_1/N+C_2\epsilon_1.$$  
 Therefore, we have 
 the diameter of $h_n\circ g_n(D_R)\leq \delta_1/N+C_2\epsilon_1$ for large $n$. 
 
 Combining with the first inequality above, we have that
 $\delta_2/2-2C_1\epsilon_2\leq \delta_1/N+C_2\epsilon_1 $.
 
  Clear if we choose  $N>\frac {\delta_1}{\delta_2/2-2C_1\epsilon_2-C_2\epsilon_1}$, 
  we get a contradiction. This expression also implies that we need to choose $\epsilon_i, i=1,2$ such that $C_1\epsilon_2+C_2\epsilon_1<<\delta_i/4, i=1, 2$.
 This finishes the proof for case (I).
 
 \vspace{2mm}
\noindent  
$\bullet $ $\bullet $ Case (II):
 Recall that in this case, $G=G_1$, $ (\Sigma, x)=({\bf C}\cup \{{\infty}\}, \infty)$. 
Each element $g\in G$ has the form $g(z)=a(z-c)$ with $a\not =0.$
Without lose generality,   we may assume that, $f_1(0)\not= f_2(\infty).$
  To justify this assumption, note that $G_1$ contains all translations. By using a translation $T$ to $f_1$, we get $f_1'=f_1\circ T$  such that  $f_1(0)\not= f_2(\infty).$ Assume that the Theorem 1.1(a) is proved for $f_1'$ and $f_2$
  with neighbourhoods $U_{\epsilon_1}(f'_1)$ and $U_{\epsilon_2}(f_2)$ and  compact subset $K'\subset G.$   Then for any $h'\in U_{\epsilon_1}(f'_1)$ and $g'\in   G\setminus K', $ we have $h'\circ g' $ is not in 
   $U_{\epsilon_2}(f_2).$

   Let  $(T^{-1})^*(U_{\epsilon_1}(f_1'))$  be the collection of all $L_k^p$-maps of the forms $\xi\circ T^{-1}$ with $\xi \in U_{\epsilon_1}(f_1').$
  Then for any $ h$ in $(T^{-1})^*(U_{\epsilon_1}(f_1'))$, we have 
  $h=h'\circ T^{-1}$ for $h'=h\circ T$ in $(U_{\epsilon_1}(f_1')).$
    Let $K=T\cdot K' $. Then for any $ g$ in $G$. Write it as $g=T\circ g'.$ we have that $g\not \in K$ if and only if $g'\not \in K'$. For any $h$ in $(T^{-1})^*(U_{\epsilon_1}(f_1'))$  and $g\not \in K$, we have that $h\circ g= h'\circ T^{-1}\circ T\circ g'=h'\circ g'$ which  is not in $U_{\epsilon_2}(f_2))$. 
  Clearly $(T^{-1})^*(U_{\epsilon_1}(f_1'))$
   can be considered as another neighbourhood of $f_1$ in ${\cal M}$. This justifies our assumption.
   
   Note that here we  have used the fact that for a fixed element $T$ in $G$, the induced action on ${\cal M}$ is a homoemorphism (actually a diffeomorphism).
  
From now on, we will assume that  $d(f_1(0),  f_2(\infty))=\delta_0>0$ and that
$\epsilon_i<< \delta_0, i=1,2.$

Now choose  a sequence of compact subsets ${\tilde K}_n\subset C^*\times {\bf C}$  defined by  ${\tilde K}_n=\{(a, c)\in C^*\times {\bf C}\,|\,\frac {1}{n} \leq |a|\leq n, \, |c|\leq n\}$. Let $K_n $ be the corresponding compact subsets in $G$ given by the map $C^*\times {\bf C}\rightarrow G$ sending $(a, c)$ to $g(z)=a(z-c)$.

Assume that the Theorem 1.1 (a) is not true so that  for all $n$ there  are  
 $g_n(z)=a_n(z-c_n)$ not in $K_n$ and $h_n$ in $U_{\epsilon_1}(f_1)$ such that  $h_n\circ
 g_n$ is in  $U_{\epsilon_2}(f_2).$

\vspace{2mm}
\noindent  
\vspace{2mm}
\noindent  
 $\bullet$ $\bullet$ $\bullet$ Claim:
 $|c_n|$ is bounded.
 
 \proof
 
 Assume that the claim is not true. Then there is a subsequence of $\{g_n\}_{n=1}^{\infty}$, denoted by the same notation, such that  $c_n\rightarrow \infty$  as $n\rightarrow \infty.$
 
Now  $h_n\circ g_n(c_n)=h_n(a_n(c_n-c_n))=h_n(0).$ Since $h_n\in U_{\epsilon_1}(f_1),$
 we have $$d(h_n\circ g_n(c_n), f_1(0))=d(h_n(0), f_1(0))$$ $$\leq \|h-f_1\|_{C^0}
 \leq C_1\|h-f_1\|_{k,  p}\leq C_1\epsilon_1.$$
 
 Since $h_n\circ g_n$ is in $ U_{\epsilon_2}(f_2),$  we have that $$d(h_n\circ g_n(c_n), f_2(c_n))\leq \|h\circ g_n-f_2\|_{C^0}
 \leq C_2 \|h\circ g_n-f_2\|_{k,  p}\leq C_2\epsilon_2.$$

 The continuity of $f_2$ at $\infty$ implies that for any given  $\epsilon$, 
 $d(f_2(c_n), f_2(\infty))\leq \epsilon$ when
  $n$ is large enough. Therefore,  $d(h_n\circ g_n(c_n), f_2(\infty))\leq  C_2\epsilon_2+ \epsilon .$       We conclude that $d(  f_1(0),   f_2(\infty))        
 \leq C_1\epsilon_1+ C_2\epsilon_2+ \epsilon=\delta_0$. This contradicts to our assumption that $\delta_0>>\epsilon_i, i=1,2.$

 \QED
 
 Therefore, we may assume that $\lim_{n\mapsto\infty}c_n=c \in {\bf C}.$
 
  This implies that 
 the sequence $\{a_n\}_{n=1}^{\infty}$ in $C^*$   either going  to zero or to infinity. As  before, we assume that $\lim_{n\mapsto\infty}a_n=\infty.$ 

 Still assume  that  diameter of the image of  $ f _2=\delta_2$. Then for  $\rho>0$  small enough,   $D^c_{\rho}(c)$, the complement  of the disc
 $D_{\rho}(c)$ of radius $\rho$ centred at $c$, 
  has the property that the diameter of $f_2(D^c_{\rho}(c))>\delta_2/2.$ 
  Moreover, there are $y_1$ and $y_2$ in $D^c_{\rho}(c)$ such that $d(f_2(y_1), f_2(y_2)))>\delta_2/2.$
  We have proved that  for large $n$, $D_{\rho}^c(c)\subset  D_{\rho/2}^c(c_n).$  
 Now let ${\tilde D}_r(\infty)$ be a small disc in $\Sigma={\bf C}\cup \{\infty\}$ centred at $\infty$ of radius $r$ measured in the Fubini-Study
 metric 
 such that the diameter of $ f_1 ({\tilde D}_r(\infty))\leq \delta_1/N.$
 Let  $D_R$ be the disc  of radius $R$ centred  at $0$ 
 such that $ {\tilde D}_r(\infty)=D^c_R$.
 
  Recall that  for any $z$ in $D_{\rho}^c(c)$, when $n$ is large enough, we have that
 $$|g_n(z)|=|a_n|\cdot |z-c_n|\geq |a_n|\cdot (|z-c|-|c_n-c|)$$ 
 $$\geq |a_n|\cdot (\rho-|c_n-c|)\geq |a_n|\cdot (\rho-\rho/2)\geq R.$$
 
 Therefore,  for large  $n$, $g_n(D_{\rho}^c(c))\subset D^c_R={\tilde D}_r(\infty) .$

 For such $n$ and $z\in D_{\rho}^c(c)$, we have that 
 
 $$d(f_1(\infty), h_n\circ g_n(z))\leq d(f_1(\infty), f_1\circ g_n(z))+d(f_1\circ g_n(z), h_n\circ g_n(z))$$ $$\leq \delta_1/N+\|f_1-h_n\|_{C^0}\leq  \delta_1/N+C_1\|f_1-h_n\|_{k, p}\leq \delta_1/N+C_1\cdot \epsilon_1.$$
 This implies that  the diameter of $h_n\circ g_n((D^c_{\rho}(c))<2(\delta_1/N+C_1\cdot \epsilon_1).$
 
  On the other hand, since $h_n\circ g_n$ is in $ U_{\epsilon_2}(f_2),$ we have that for $y_1$ and $y_2$ in $D^c_{\rho}(c)$ chosen above, $$ d(h_n\circ g_n(y_1), h_n\circ g_n(y_2))$$ $$ >d(f_2(y_1), f_2(y_2)))-d(f_2(y_1), h_n\circ g_n(y_1))-d(h_n\circ g_n(y_2), f_2(y_2))     $$ 
  $$ > \delta_2/2-2C_2\cdot \epsilon_2.$$ This implies that the diameter of $h_n\circ g_n ((D^c_{\rho}(c))>\delta_2-2C_2\cdot \epsilon_2.$

 We conclude that $\delta_2-2C_2\cdot \epsilon_2<$ the diameter of $h_n\circ g_n((D^c_{\rho}(c))<2(\delta_1/N+C_1\cdot \epsilon_1).$
   By  choosing  $N>>\frac {2\delta_1}{\delta_2-2C_2\cdot \epsilon_2-2C_1\cdot \epsilon_1}, $ we get a contradiction. Again we should
 choose $\epsilon_i << \delta_j, i, j=1, 2.$

This finishes the proof of the Theorem 1.1 and 1.2.
\QED

\begin{remark}
\vspace{2mm}
\noindent  
\vspace{2mm}
\noindent
(A) There are  corresponding statements for all Theorems and corollaries 
for $\Sigma=S^n$ with its standard metric and $G$ to be the conformal group of
of $S^n$. The proof for these general theorems will be given in a subsequent paper. In fact, the proof goes along  similar line as we did here.

\vspace{2mm}
\noindent  
\vspace{2mm}
\noindent 
(B) In the proof above, we have used diameter function, $ diameter (h(\Sigma'))$
 with $(h, \Sigma')\in {\cal M}\times {\cal S}(\Sigma), $ where ${\cal S}(\Sigma)$ is the collection of all sub-surfaces of $\Sigma.$  There are many other functions, such as energy function  can be used to give essentially the same proof.

\end{remark}

The original  motivation of  this  work is to prove the much weaker but more reasonable result that if a $L_k^p$-map $f$ is stable in the sense that (i)either $\omega (f)>0$, or (ii) it  has no  "infinitesimal" automorphism, then $f$ is
 $G$-stable in the sense described in Theorem 1.1. In other word, we want to find
  the corresponding statement in our infinite dimensional setting of the  
  well-known result by  Gieseker and  Mumford  that 
  for stable curves  stability in the sense of the Deligne-Mumford is equivalent to the GIT
  stability.   It turns  out that the proof for this weaker result also works only under the assumption that $f$ is not constant, which leads to the Theorem 1.1.

\section{ More Corollaries  of the Theorem 1.1}

\begin{cor} 
Given     
 two maps $f_1$ and $f_2$ in ${ {\cal M}}^*$  not in the same $G$-orbit,  there exist small G-neighbourhoods $U^G_{\epsilon_1}(f_1)=G\cdot U_{\epsilon_1}(f_1)$ and $U^G_{\epsilon_2}(f_2)=G\cdot U_{\epsilon_2}(f_2)$  such that they  do not intersect each other.   In other words, the $G$-space  ${ {\cal M}}^*$ is $G$-Hausdorff.

\end{cor}

 \vspace{2mm}
\noindent 
 \proof
 
 \vspace{2mm}
\noindent
 By Theorem 1.1, for any $g\not \in$ the compact set $ K_1$ and $h\in  U_{\epsilon_1}(f_1)$, $h\circ g$ is not in $U_{\epsilon_2}(f_2)$. By our assumption, we may assume that 
 $  U_{\epsilon_1}(f_1)$ and $U_{\epsilon_2}(f_2)$ have no intersection.
 
 \vspace{2mm}
\noindent
\vspace{2mm}
\noindent $\bullet$ Claim: when $\epsilon_i,i=1,2$  are small enough, 
 $U^G_{\epsilon_1}(f_1)\cap U_{\epsilon_2}(f_2)$ is empty.

\proof

\vspace{2mm}
\noindent
If this is not true, there  are $h_i\in U_{\delta_i}(f_1)$ and $g_i\in K_1$
such that $h_i\circ g_i$ is in $ U_{\delta_i}(f_2)$ with $\delta_i\mapsto 0.$
 The compactness of $K_1$ implies that after taking a subsequence, we have that
 $\lim_{i\mapsto \infty}g_i=g\in K_1.$ Since $\delta_i\mapsto 0$, we have that
 $f_1=\lim_{i\mapsto \infty}h_i$ and $f_2=\lim_{i\mapsto \infty}h_i\circ g_i=f_1\circ g.$
Hence, $f_1$ and $f_2$ are in the same orbit which contradicts to our assumption.
Note that in the last identity above, we have used the fact that the action map
$\Psi: G\times {\cal M}\rightarrow {\cal M}$ is continuous.

\QED

Of course the same proof also implies that 
$U^G_{\epsilon_2}(f_2)\cap U_{\epsilon_1}(f_1)$ is also empty for sufficiently
 small $\epsilon_i, i=1, 2.$

If $h\in U^G_{\epsilon_1}(f_1)\cap U^G_{\epsilon_2}(f_2), $  then there are
$h_i\in U_{\epsilon_i}(f_i) $ and $g_i\in G, i=1, 2$ such that 
$h=h_1\circ g_1=h_2\circ g_2. $ Hence $h_2=h_1\circ g_1\circ g_2^{-1}$ and
$U^G_{\epsilon_1}(f_1)\cap U_{\epsilon_2}(f_2)$ is not empty. This contradicts
to the above claim.

\QED

This implies  the following corollary.
  
 \begin{cor} 
The $G$-space  ${ {\cal M}}$ is $G$-Hausdorff.

\end{cor}
 
\vspace{2mm}
\noindent

\proof 

  We have already  proved that if $f_1$ is a constant map and $f_2$ is not, the above corollary is still true for  much  simpler reason. 
  Therefore we only need to show that if $f_1\not =f_2$ are two constant maps, then above  corollary is still true.
  
  To this end, let $B_{\epsilon'_1}(c_1)$ and  $B_{\epsilon'_2}(c_2)$ be two open balls in $M$, which do not intersect. Here $c_1$ and $c_2$ are the values of the two constant maps $f_1$ and $f_2$ respectively. Clearly if $||h_i-f_i||_{C^0}=max_{x\in \Sigma}|h_i(x)-c_i|<\epsilon'_i, i=1, 2, $ then the image of $h_i$ is contained in $B_{\epsilon'_i}(c_i).$ Moreover, since for  any  $h_i$ and $g_i\in G$, 
   the  image of $h_i\circ g_i=$ the  image of $h_i,$ for any  $h_1$ and $h_2$ as above, their $G$-orbits
  $G\cdot h_1$ and  $G\cdot h_2$ do not intersect. Clearly by our assumption   for $\epsilon_i<<
  \epsilon'_i,$
   any  $h_i, i=1, 2$ in $U_{\epsilon_i}(f_i)$ satisfies the condition $||h_i-f_i||_{C^0}<\epsilon'_i$, hence $G\cdot U_{\epsilon_1}(f_1)$
  and $G\cdot U_{\epsilon_2}(f_2)$ do not intersect.
  
  QED

\begin{cor} 
Given any    
 $f$ in ${ {\cal M}}$, the G-orbit $G\cdot f$ is closed in ${\cal M}$.  
\end{cor}
 
\proof

We only need to consider the case that $f$ is a non-constant map.

The proof is similar to the proof of  the  first corollary. 
Rename $f$ as  $f_1$. If the corollary is not true, there exist $g_i\in G$ and $f_2\in {\cal M}^*$ such that
$f_2=\lim_{i\mapsto \infty}f_1\circ g_i,$ but $f_2$ is not in $G\cdot f_1$.
Therefore for any $U_{\epsilon_2}(f_2), $  when $i$ is large enough,  $f_1\circ g_i$ is in $U_{\epsilon_2}(f_2). $
On the other hand,  the Theorem 1.1 with the same notation there implies that for all such $i$,  $g_i$ is in the compact set $K_1$. Therefore, we may assume that
$\lim_{i\mapsto \infty} g_i=g$ in $K_1.$ Consequently, $f_2=\lim_{i\mapsto \infty}f_1\circ g_i=f_1\circ g.$ That is $f_2\in G\cdot f_1$ which is  a contradiction.

\QED

 Essentially the same argument proves the following stronger result.
 \begin{cor} 
Given any  non-constant map   
 $f$ in ${ {\cal M}}$, there is a small closed ${\delta}$-neighbourhood $B_ {\delta}(f)$ such that the G-orbit $G\cdot B_ {\delta}(f)$ is closed in ${\cal M}^*$.  In other words, ${ {\cal M}}^*$ is $G$-regular in the sense that for any $G$-closed subset $C$ in ${ {\cal M}}^*$  and $f\not \in C$ , there are $G$-open neighbourhoods $U_1$ and $U_2$ of $C$ and $G\cdot f$ respectively such that $U_1$ and $U_2$ do not intersect.

\end{cor} 
\proof

Rename $f$ as  $f_1$.  If the corollary is not true, then for some $\delta>0, $ there exist  sequences $\{h_{ j}\}_{j=1}^{ \infty}\in  B_ {\delta}(f_1)$ and  $g_{ j}\in G$ and $f_{2}\in {\cal M}^*$ such that
$f_{ 2}=\lim_{j\mapsto \infty}h_{ j}\circ g_{j},$ but $f_{ 2}$ is not in $G\cdot B_ {\delta}(f_1)$.

   For  any $U_{\epsilon_2}(f_{ 2}), $  when $j$ is large enough,  $h_j\circ g_{ j}$ is in $U_{\epsilon_2}(f_{2}). $
On the other hand,  the Theorem 1.1  implies that  for proper choices   of   the radius ${\delta}$ of  $B_ {\delta}(f_1)$ and the radius ${\epsilon_2}$ of $U_{\epsilon_2}(f_{ 2}) $,    $g_{ j}$ is in the compact set $K_1$ defined in the Theorem 1.1. Therefore, we may assume that
$\lim_{i\mapsto \infty} g_{ j}=g$ in $K_1.$  Hence
$$\lim_{j\mapsto \infty} h_{ j}=\lim_{j\mapsto \infty} (h_{ j}\circ g_{j})\circ  \lim_{j\mapsto \infty} g_{j}^{-1} $$ $$ =f_{ 2}\circ  g^{-1}.$$
 Since $ B_ {\delta_i}(f_1)$  is closed, we have $\lim_{j\mapsto \infty} h_{ j},$ denoted by $h$,  is in  $ B_ {\delta}(f_1)$.
Therefore, $f_{ 2}=\lim_{j\mapsto \infty}h_{ j}\circ\lim_{j\mapsto \infty} g_{j}=h\circ g\in  G\cdot B_ {\delta}(f_1)$. This is a contradiction.

\QED

   A  corollary of the above results is the following.
  \begin{cor} 
   
 The quotient space  ${{\cal M}}/G$ of the unparametrized
 $L_k^p$-maps is Hausdorff. In particular,
 the space   ${\cal B}={\tilde {\cal B}}/G$  of the unparametrized stable $L^p_k$-maps  with fixed domain  is Hausdorff.
  Here a  $L_k^p$-map $f:\Sigma\rightarrow M$ with $\Sigma$ being  a  node curve is said to be stable if (i)any of its genus zero "free" component is non-trivial; (ii) the  stabilizer   $Stab_f$  is  finite. 
  \end{cor}

In the sequel of this paper [L1], we will generalize the above result to the general stable $L^p_k$-maps allowing varying
of the domain $\Sigma$ in the Deligne-Mumford type of moduli spaces. The results in this section and their  generalizations will be used
for the regularization of the moduli space of stable $J$-holomorphic  curves described in [L].

A weaker statement that the genus zero moduli space ${\cal M}(J, A)$ of unparametrized stable  $J$-holomorphic maps of class $A\in H_2(M, {\bf Z})$ is Hausdorff was proved in [LT]. The proof there also works for the higher genus case.
The Hausdorffness  for  a neighbourhood  in ${ {\cal B}}$ of  unparametrized stable  $L_k^p$-maps covering  the space ${\cal M}(J, A)$  was used implicitly in [LT] without proof.
To author's knowledge,   the  Hausdorffness  for   ${ {\cal B}}$ of  unparametrized stable  $L_k^p$-maps was first proved  by  Hofer,  Wysocki and Zehnder  in [HWZ] under  a stronger   notion  of  stability. In particular,  to define the stability in the sense of [HWZ],   the target space $M$  is required to be  a symplectic manifold. We refer the readers to  [HWZ] for the definition.

\vspace{2mm}
\noindent 
\vspace{2mm}
\noindent  
Note: There is a parallel discussion to the results so far  for the case that $G_0=SL(2, {\bf R})$
 acting on the  $\Sigma\simeq {\bf H}$= upper half-plane $\simeq D^2=$ closed disc as  well as the case for    its subgroup $G_i, \, i=1,2$
 preserving one or two marked points on the boundary. This together with its generalization to the case that the domain $\Sigma$ has deformation appeared in Lagrangian Floer homology will be treated in [L2].

\vspace{2mm}
\noindent 
\vspace{2mm}
\noindent  
There are a few immediate related  results and questions that will be   treated in a subsequent paper. 

\vspace{2mm}
\noindent  
(A)  Let ${\cal M}^{**}$  be the subspace of ${\cal M}^{*}$ consisting of all $L_k^p$-maps  $f$ whose stabilizer $Stab_f$ is finite.
Then by using  the  argument   in [La], one can show   that the  quotient  space ${\cal M}^{**}/G$   is paracompact. This implies  that ${\cal M}^{**}$ is $G$-paracompact.

\vspace{2mm}
\noindent  
(B) According to the homology classes  represented by its elements, the spaces ${\cal M},$ ${\cal M}^{*}$ and ${\cal M}^{**}$ are decomposed as:
${\cal M}=\cup _{A\in H_2(M, {\bf Z})}{\cal M} (A)$, etc. Note  that for any $A\not = 0, $ ${\cal M} (A)={\cal M}^* (A).$
In this case, ${\cal M} (A)={\cal M}^* (A)$ is decomposed further into two essential strata ${\cal M}^{i} (A), i=0, 1,$ according to   the  dimension of $Stab_f=0 $ or $1$.  The reason for this is that by the corollary 1.1,  in this  case, $Stab_f$ is a compact Lie subgroup of $G$. On the hand, upto the conjugations, the only connected  compact subgroup of $SL(2,{\bf C})$ are $SU(2)$ or $S^1$. Since $SU(2)$ acts on $\Sigma \simeq S^2$ transitively, the condition that  $SU(2) $ is in $ Stab_f$ implies that $f$ is a constant map. Therefore, for any  non-constant map $f$,  the connected component of identity of  $Stab_f$ is either    trivial  or equal to $S^1$ upto a conjugation.

\vspace{2mm}
\noindent  
(C)  The simplest case for   $f$ in ${\cal M}^{1} (A)$ is that $f:S^2\rightarrow M$ factors through as $f={\bar f}\circ \pi. $ Here $\pi:S^2\rightarrow [-1, 1]$ is the quotient map defined by sending each closed orbit of the "standard" $S^1$-action  on $S^2(1)\subset {\bf R}^3$ to the intersection of the rotation axis with the plane passing through the orbit, and ${\bar f}:[-1, 1]\rightarrow M$ is the obvious induced map. In the case that $stab_f$ is connected, upto an conjugation from an element $G$, this simplest case is the "normal" form for the case here. Of course, for the general case,  the corresponding  $\pi$ is more complicated as we allow a further   finite  equivalence  relation.  However, the possibility to have "normal" form in the above simplest case suggests that the question here is "discrete", hence it makes sense  to try  to find the "normal" forms  even for the general case. 
A closely related  question is  to classify  the singular line foliation $\xi_f$ on $S^2$  upto a proper equivalence relation (such as conjugation by a differomorphism)  whose generic integral curves are all closed $S^1$. Here the foliation  $\xi_f$ is defined to be  the  kernel of $df.$

\vspace{2mm}
\noindent  
(D) For each stratum above, in particular  for the stratum ${\cal M}^{1} (A), $ the usual $G$-equivariant   tubular neighbourhood theorem
stated, for instance  in [GGK] page 180, can be generalized  to  this case accordingly by using weakly smoothness or $sc$-smoothness. The resulting $G$-equivariant   tubular neighbourhood of a $G$-oribt $G\cdot f$ in ${\cal M}^{1} (A) $ is essentially an principal $S^1$-bundle ( upto a finite equivalence relation)   over the orbit in the proper category.

\vspace{2mm}
\noindent  
(E)  The result in (C) can be used to show that the quotient space ${\cal M}^{*}/G$ is paracompact  so that
${\cal M}^{*}$   is $G$-paracompact.

\vspace{2mm}
\noindent  

\vspace{2mm}
\noindent 
$\bullet$ $\bullet$ On the notion of proper $G$-action:

\vspace{2mm}
\noindent 
The reader might have already noticed that the theorems and corollaries above are the analogies of the familiar ones in the case  that $G$ acts properly  on a finite dimensional manifold.  We close this section  with a  comparison of the  definition  of properness of $G$-action defined here with the one used in finite dimensional situation.

To this end, let $M$ be  a finite dimensional manifold acted by a a non-compact Lie group $G$. Recall that the action of  $G$ is said to be proper if the total action map $\Psi: M\times G\rightarrow M\times M$  given by $\Psi(m ,g)=(m, m \cdot g)$ is proper.

We now make two reductions. 

\vspace{2mm}
\noindent 
(I) To check the properness,  it is sufficient to look  at all compact set $K$ in $M\times M$ of the form $K=K_1\times K_2$ with $K_i, i=1, 2 $ is contained in $M$. To justify this, note that   for any compact set $K\in M\times M$, $\pi_1(K)\times \pi_2(K)$ is compact in $M\times M$,  hence $\Psi ^{-1}(\pi_1(K)\times \pi_2(K))$ is   compact in $M\times G$. 
Here $\pi_i$ is the projections of $M\times M$ to its two factors.
Since  $\Psi ^{-1}(K)$ is  a closed subset in $\Psi ^{-1}(\pi_1(K)\times \pi_2(K))$, it is compact. 

\vspace{2mm}
\noindent 
(II) For any compact subset $K_1\times K_2$ in $M\times M$, $\Psi ^{-1}(K_1\times K_2)$ is compact if and only if its projection to $G$, $\pi_G(\Psi ^{-1}(K_1\times K_2))$ is compact in $G$.

\vspace{2mm}
\noindent 
Therefore in the finite dimensional case, the $G$-action is  proper if and only if for any compact subset $K_1\times K_2$ in $M\times M$, 
$\pi_G(\Psi ^{-1}(K_1\times K_2))$ is compact in $G$.

\vspace{2mm}
\noindent 
In our infinite dimensional setting, the  Theorem 1.1 can be reformulated as follows.

\vspace{2mm}
\noindent 
$\bullet$ $\bullet$ The equivalent forms of the Theorem 1.1.

\vspace{2mm}
\noindent 
(A) For any $(f_1, f_2) \in {\cal M}^*\times {\cal M}^*, $ there exists a product  open neighbourhood  $U_{\epsilon}=U_{\epsilon_1}(f_1)\times U_{\epsilon_2}(f_2)$  such that 
$\pi_G(\Psi ^{-1}(U_{\epsilon}))$ is pre-compact in $G$.

\vspace{2mm}
\noindent 
 A standard argument in this situation implies  a stronger form of the above statement.

\vspace{2mm}
\noindent 
(B)  For any compact subset $K = K_1\times K_2 \in {\cal M}^*\times {\cal M}^*, $ there exists a product  open neighbourhood  $U_{\epsilon}=U_{\epsilon_1}(f_1)\times U_{\epsilon_2}(f_2)$  such that 
$\pi_G(\Psi ^{-1}(U_{\epsilon}))$ is pre-compact in $G$.

\vspace{2mm}
\noindent  Of course,  in  (B) above, one can simply use $K$ and $U_{\epsilon}$  not necessarily to be   a  product of two sets.
However, the proof in [L0] seems suggesting that, at least in the case that $  {\cal M}^* $  is the space of $v$-stable maps (defined in next section), $U_{\epsilon_1}(f_1)$ can be an arbitrary "bounded" open neighbourhood. Whether or this is true   will be  decided in a subsequent paper.

\section{Generalizations of the Main Theorem}

In this section, we will state two theorems that generalize  the  Theorem 1.1.   All the other results in this paper have the corresponding ones in this general setting. 

We  assume that $G$ is a connected reductive group. Let ${ {\cal M}}^*$ be the collection of all $L_k^p$-maps $f:\Sigma\rightarrow M$ which are $v$-stable.
Here we assume that $m_0=k-\frac{n}{p}>1, $ where $n=\dim (\Sigma).$

\begin{definition}

A $L_k^p$-map $f:M\rightarrow N$ is said be to $v$-stable if its volume $v(f)>0$.

\end{definition}

For the first case, $\Sigma$ is a smooth projective manifold, and $G$ acts on $\Sigma$ algebraically. Then we have the following theorem.

\begin{theorem}

 Assume that  $G$  is a connected  reductive group   acting on a smooth projective variety  $\Sigma $ in the  way described as above. 
 Then  the  action  of $G$  on  ${ {\cal M}}^*$ is  proper in the following sense: 
  Given any two maps $f_1$ and $f_2$ in ${ {\cal M}}^*$, there exist small neighbourhoods $U_{\epsilon_1}(f_1)$ and $U_{\epsilon_2}(f_2)$ of $f_1$ and $f_2$ in ${ {\cal M}}^*$ and two compact sets $K_{f_1}$ and $K_{f_2}$ in $G$ such that
  
  \vspace{2mm}
\noindent
  (a) for any $h$ in $U_{\epsilon_1}(f_1)$ and $g$ in $G\setminus K_{f_1}$, 
  we have that $h\circ g$ is not in $U_{\epsilon_2}(f_2)$;
  
   \vspace{2mm}
\noindent
(b) the corresponding statement holds for any $h$ in $U_{\epsilon_2}(f_2)$ and $g$ in $G\setminus K_{f_2}$.
  
  \end{theorem}

 Clearly by our assumption, we may assume  that  $\Sigma$  is embedded in some $ {\bf CP}^k$  not lying in any of its  hyperplane sections. 
The key  step to prove this theorem is   to note the  well-known fact that in this situation   for any element $g\in G$, the action of  $g$ on $\Sigma$ is given by a element in $PSL(k+1, {\bf C})$. This linearises the situation. The rest of the proof is elementary using similar ideas as this paper. But instead of using diameter function in this paper, for the proofs  of  this and next theorems,  the volume function is used as  expected from the definition of $v$-stability.

It is possible to give a  proof for the above theorem without using the  complete information of the  ambient  space $ {\bf CP}^k$, but only  assuming that the action of $G$ on $\Sigma$  comes from a moment map in the sense that the induced action of  the fixed maximum compact subgroup $K$ is generated by a moment map. This proof relies on the   works  of Atiyah in [A] and Giullemin and Sternberg in [GS] on the images of the moment maps of tours actions and  the work of Kirwan in [K]. It leads to the following generalization.

\begin{theorem}

 Assume that  $G$  is a connected  reductive group   acting on a compact and smooth  Kahler manifold   $\Sigma $ holomorphcally such that the action comes from a moment map.
 Then  the  action  of $G$  on   ${ {\cal M}}^*$ is  proper in the  sense described in the above theorem.
  \end{theorem}

At the writing of this paper, only the proof for the first theorem above is completely carried out. The author does not expect real difficulties  for the proof of the second theorem.

\begin{remark}

 These theorems are supposed to be the generalizations for the corresponding statements in Theorem 1.1.  However, in  Theorem 1.1  , in stead of assuming that each element  $f$ is $v$-stable, which means that  the area of $f$ is larger than zero, we only assume that $f$ is a non-constant map.
 Of course this last assumption is equivalent to the assumption that the energy  
of each  $f$ is larger that than zero,  which  is weaker than the condition that $f$ is $v$-stable.

Therefore, one may ask if the theorems above in this paper  are still true under the  weaker  condition that each $f$ is  a non-constant  map. Equivalently,  we  can   define the $m$-dimensional energy $v_1(f)=\int_{\Sigma} \|df\|^m d\nu_{\Sigma}$ as well as 
the corresponding notion of $f$  being a $v_1$-stable map.  Here $\|df\|(x)$ is 
   the norm of the linear map $df(x):T_x\Sigma \rightarrow T_{f(x)}M$ measured by  the metrics on $T_x\Sigma$ and $T_{f(x)}M.$   Then  we want to know if  the space of $v_1$-stable maps is $G$-stable in the sense of the above theorems.
   It turns out that this straightforward generalization using $v_1$-stability is  not the right one. However, as the proof in [L0]shows,   the right  condition for the generalization is to require that  $f$ is not constant with respect to the fixed point sets of all subgroups $S^1$ (or their complexificatons $C^*$ ) in $G$. More specifically, let $S\subset \Sigma$ be a fixed point set of a subgroup $S^1$ in $G$,
   then the requirement is that the image of $f:\Sigma\rightarrow M$ is not contained in the image of 
   $f|_S$ for all such $S$. Clearly in the case that $\Sigma={\bf CP}^1$ and 
   $G=SL(2, {\bf C})$,  any  fixed point set $S$  always consists of two points in $\Sigma$ so that in this case the above  condition on $f$ is equivalent to the condition   in  Theorem 1.1.  However, for the general case here, above optimal condition is not easy to check, hence may not be very useful.
   
\end{remark}

\section{ A Pseudo-moment Map for the Action  of  $K$ on ${\cal M}^*$}
 
  In the following we only give the definitions of the pseudo-moment maps without  detailed explanations. The details  will be given in the revised version of this paper.

 We start with the general case described in the second theorem of Sec. 4. Fix a maximum compact subgroup $K$ in the reductive group $G$ acting holomorphically on the smooth and compact Kahler manifold $\Sigma.$ Denote the Lie algebra of $K$ by ${\bf k}$. Assume that $\mu:\Sigma\rightarrow {\bf k}^*$ be the moment map that generates the $K$-action. We now define the induced pseudo-moment map ${\bf m}:{\cal M}^*\rightarrow {\bf k}^*$
 as follows. For any  $\xi$ in ${\bf k}$, let $c_{\xi}=(max (\mu_{\xi})-min (\mu_{\xi}))/2$ be the average of the maximum and minimum of the function $\mu_{\xi}=<\mu, \xi>$ defined on  $\Sigma $. Denote the inverse image in $\Sigma$, $\mu^{-1}((-\infty, c_{\xi}])$ by $\Sigma^{-}_{\xi}$. Assume that $c_{\xi}$ is a regular value first so that  $\Sigma^{-}_{\xi}$ is a compact sub-manifold with boundary (a closed region) in $\Sigma$.

For any $f:\Sigma\rightarrow M$ in ${\cal M}^*$ and $\xi$ in ${\bf k}$, we define the pseudo-moment map ${\bf m}:{\cal M}^*\rightarrow {\bf k}^*$ by the identity:  $<{\bf m}(f), \xi>=v(f)/2-v(f|_{\Sigma^{-}_{\xi}}).$ Here as before, $v(f)$ and $v(f|_{\Sigma^{-}_{\xi}})$ are the volumes of the two  maps. In the case that $c_{\xi}$ is not a regular value, by taking a sequence of regular values approaching to $c_{\xi}$
and applying the above definition,  we define $<{\bf m}(f), \xi>$  to be  the limit of the resulting sequence provided that we can prove that the 
limit exists  independent of the choices made above.

\vspace{2mm}
\noindent 
\vspace{2mm}
\noindent  Note: There is a minor inaccuracy in above definition.   The  range  of ${\bf m}$ is ${\hat {\bf k}}^*$ rather than ${\bf k}^*$. Here ${\hat {\bf k}}$ is the oriented real blow-up of ${ {\bf k}}$ at the origin and hence is a "ray" bundle over a sphere, and  ${\hat {\bf k}}^*$ is its fiber-wise "dual".

In the case that $\Sigma= S^2$, $G=PSL (2)$ and $K=PSU(2)=SO(3).$ There is a slight different and more geometric description of above pseudo-moment map.  We  can  work  at the group level rather than using Lie algebra ${\bf so}(3)$. Consider   the  double covering   $SL (2)$ of $PSL (2)$ and $SU(2)$  of $SO(3).$ For any $g\not=\pm$ identity in $SU(2)$ its image in $SO(3)$ is a non-trivial rotation of $\Sigma=S^2$ sitting inside ${\bf R}^3.$ The element $g$ picks up an orientation for the   axis of the rotation. Therefore it makes sense to talk about the lower-half
sphere $\Sigma^{-}_{g}$ with respect this oriented axis, which corresponds to  $\Sigma^{-}_{\xi}$ in the above general description.

The   pseudo-moment map in this description is  ${\bf m}:{\cal M}^*\times {\widetilde SU}(2)  \rightarrow {\bf R}$  defined by ${\bf m}(f, g)=v(f)/2-v(f|_{\Sigma^{-}_{g}}).$
Here 
${\widetilde SU}(2)$ is the oriented real blow-up of $SU(2)\simeq S^3$ at $\pm$identity. 

As  mentioned in the introduction, the motivation to introduce these  pseudo-moment maps is to  get a global slice for the non-compact  directions of the $G$ action on a neighbourhood in ${\cal M}^*$ of the moduli space of stable $J$-holomorpfic maps. The details  on how this can be done as  well as  the applications  of the existence of such a global slice to GW-theory  will be discussed in the  revised version of this paper.

\end{large}
\end{document}